\documentclass[12pt]{article}
\usepackage{amsmath,amsthm,amsfonts}
\usepackage{fullpage}

\newtheorem{theorem}{Theorem}

\title{On the context-freeness of the set of words containing overlaps}
\author{Narad Rampersad \\
David R. Cheriton School of Computer Science \\
University of Waterloo \\
Waterloo, Ontario N2L 3G1 (Canada) \\
{\tt nrampersad@math.uwaterloo.ca}
}

\begin{document}
\date{\today}
\maketitle

\begin{abstract}
We show that the set of binary words containing overlaps is not
unambiguously context-free and that the set of ternary words containing
overlaps is not context-free.  We also show that the set of binary words that
are not subwords of the Thue--Morse word is not unambiguously context-free.
\end{abstract}

\section{Introduction}

An \emph{overlap} is a word of the form $axaxa$, where $a$ is a single letter
and $x$ is a (possibly empty) word.  A word is \emph{overlap-free} if
it does not contain an overlap as a subword.  Let $L_o(k)$
denote the language of all words over the alphabet $\{0,1,\ldots,k-1\}$
that contain an overlap as a subword.

By applying the interchange lemma \cite{ORW85,RW82}, Gabarr\'o
\cite{Gab85} proved that for $k \geq 4$, $L_o(k)$ is
not context-free, thus partially solving an open problem of Berstel
\cite{Ber84}.  Since the work of Gabarr\'o, it has remained
an open problem to determine whether
or not $L_o(2)$ and $L_o(3)$ are context-free.
We show that $L_o(2)$ is not unambiguously context-free and that
$L_o(3)$ is not context-free.  We also show that the set of binary words that
are not subwords of the Thue--Morse word is not unambiguously context-free.

\section{Overlap-free words}
\label{overlap-free}

In this section we review some standard results concerning binary
overlap-free words.

Let $\mu$ denote the \emph{Thue-Morse morphism}, that is, the
morphism that maps $0 \to 01$ and $1 \to 10$.  The \emph{Thue-Morse word}
\[
\mathbf{t} = \mu^\omega(0) = 0110100110010110\cdots
\]
is well-known to be overlap-free \cite{MH44,Thu12}.

Let \[A = \{00,11,010010,101101\}.\]
Pansiot \cite{Pan81} and Brlek \cite{Brl89} proved the set of squares
in $\mathbf{t}$ is exactly \[\mathcal{A} = \bigcup_{k \geq 0}\mu^k(A).\]
Using this result, one easily shows (see, for example, \cite{CRS05}) that
for any position $i$, there is at most one square in $\mathbf{t}$
beginning at position $i$.

\section{Binary words containing overlaps}

\begin{theorem}
\label{ambig}
The language $L_o(2)$ is not unambiguously context-free.
\end{theorem}

\begin{proof}
We will need the following result due to Fatou \cite{Fat06}
(a more convenient reference may be \cite[Part~VIII, Chap.~3, No.~167]{PS71};
a stronger result was conjectured by P\'olya and proved by Carlson
\cite{Car21}): \emph{A power series $\sum_{n \geq 0} a_n z^n$ with integer
coefficients and radius of convergence $1$ is either rational
or transcendental over $\mathbb{Q}(X)$.}

Let
\[
F(X) = \sum_{n \geq 0} a_n X^n
\]
be the generating series of the overlap-free words.  That is,
$a_n$ is the number of overlap-free words of length $n$
over a two letter alphabet.  By the Chomsky--Sch\"utzenberger
Theorem \cite{CS63} (see \cite[Chap. 16]{KS86} for a proof; see also,
for example, \cite{All99,Fla87} for applications
to other languages), if $L_o(2)$ is unambiguously
context-free, then $F(X)$ is algebraic over $\mathbb{Q}(X)$.
To prove the theorem it suffices then to show that $F(X)$ is
transcendental.

We will need the following result due to Lepist\"o \cite{Lep96} on the
enumeration of overlap-free words (compare also the earlier
work of Restivo and Salemi \cite{RS85}, Kfoury \cite{Kfo88},
Kobayashi \cite{Kob88}, and Cassaigne \cite{Cas93}):
\begin{equation}
\label{enum}
a_n = \Omega(n^{1.217}) \quad\text{and}\quad a_n = O(n^{1.369}).
\end{equation}

Since $a_n = O(n^{1.369})$, $F(X)$, as a complex power series, has radius
of convergence $1$, and so by Fatou's theorem is either rational
or transcendental over $\mathbb{Q}(X)$.  To complete the proof we
must show that $F(X)$ is not rational.  If $F(X)$ were rational, then
the coefficients $a_n$ could be written in the form
\[
a_n = \sum_{i=1}^m A_i(n) \alpha_i^n,
\]
for some $m$, where $\alpha_i$ is a characteristic root of
multiplicity $n_i$ of the linear recurrence satisfied by
$(a_n)_{n \geq 0}$, and $A_i(X)$ is a polynomial of degree at most $n_i - 1$
(see \cite[Section~1.1.6]{EPSW03}).
But from (\ref{enum}) we see that this is not possible, so $F(X)$ is not
rational and the proof is complete.
\end{proof}

We conclude this section by considering a variation on the language $L_o(k)$.
Given a word $w$, and writing $w = xy$, we say the word $yx$ is a
\emph{conjugate} of $w$.  Let $\widetilde{L_o}(k)$
denote the language of all words $w$ over the alphabet $\{0,1,\ldots,k-1\}$
such that some conjugate of $w$ contains an overlap as a subword.

\begin{theorem}
The language $\widetilde{L_o}(2)$ is not unambiguously context-free.
\end{theorem}

\begin{proof}
Harju \cite{Har86} showed that the binary circular overlap-free words
have lengths of the form $2^n$ or $3 \cdot 2^n$, $n \geq 0$.
The generating series of the complement of $\widetilde{L_o}(2)$ is a thus
a so-called ``gap series'' (or ``lacunary series'').  By Hadamard's
gap theorem \cite[Theorem~16.6]{Rud86} it admits its circle of convergence
as a natural boundary and hence is transcendental.  Applying
the Chomsky--Sch\"utzenberger Theorem, we conclude that $\widetilde{L_o}(2)$
is not unambiguously context-free.
\end{proof}

\section{Ternary words containing overlaps}

In this section we adapt the argument of Gabarr\'o \cite{Gab85} to
prove the following theorem.

\begin{theorem}
\label{ternary}
The language $L_o(3)$ is not context-free.
\end{theorem}

Before beginning the proof, we recall the interchange lemma \cite{ORW85}.

\begin{theorem}[Ogden, Ross, and Winklmann]
Let $L \subseteq \Sigma^*$ be a context-free language.  There exists a
constant $c$, depending only on $L$, such that for all $n \geq 2$, all
subsets $R \subseteq L \cap \Sigma^n$, and all $m$, $2 \leq m \leq n$,
there exists a subset $Z \subseteq R$, $Z = \{z_1,z_2,\ldots,z_k\}$,
such that
\begin{itemize}
\item[(a)] $k \geq \frac{|R|}{c(n+1)^2}$;
\item[(b)] $z_i = w_i x_i y_i$, $1 \leq i \leq k$;
\item[(c)] $|w_1| = |w_2| = \cdots = |w_k|$;
\item[(d)] $|y_1| = |y_2| = \cdots = |y_k|$;
\item[(e)] $m/2 \leq |x_1| = |x_2| = \cdots = |x_k| \leq m$;
\item[(f)] $w_i x_j y_i \in L$, $1 \leq i,j \leq k$.
\end{itemize}
\end{theorem}

\begin{proof}[Proof of Theorem~\ref{ternary}]
Let $n = 2^{2k+1}+1$ for some $k \geq 0$.  Let $x = \mu^{2k}(0)$ and
let $w = 0xx$.  Then $w$ is an overlap, but no proper subword of $w$
is an overlap.  To see this, note that $xx$ is a subword of the Thue-Morse
word and is therefore overlap-free.  Any overlap contained in $w$
must therefore begin from the first position of $w$.  If $w$ begins
with two distinct overlaps, then $xx$ begins with two
distinct squares, contradicting the observation made in
Section~\ref{overlap-free}.

Suppose that $L_o(3)$ is context-free.
Let $\psi$ be the morphism defined by $\psi(0) = 0$ and
$\psi(1) = \psi(2) = 1$.  Define
\[
R = \{ 0yy : y \in \psi^{-1}(x) \}.
\]
Note that $|R| = 2^{(n-1)/4}$.  Applying the interchange lemma,
we see that there exists $Z \subseteq R$ with
\begin{equation}
\label{upperbnd}
|Z| \geq \frac{2^{(n-1)/4}}{c(n+1)^2}.
\end{equation}

Choosing $m = (n-1)/2$, and recalling that if $z_i = w_ix_iy_i \in Z$,
then $m/2 \leq  |x_i| \leq m$, we see that $w_ix_jy_i \in L$ only if
$x_i = x_j$.  Fixing $x_i$, we easily verify that there are at most
$2^{(n-1)/8}$ words $w_jx_jy_j$ with $x_i = x_j$, so that
$|Z| \leq 2^{(n-1)/8}$, contradicting (\ref{upperbnd}) for $n$
sufficiently large.  This concludes the proof.
\end{proof}

\section{Generalized Thue--Morse words}

In this section we show that the set of binary words that are not
subwords of the Thue--Morse word $\mathbf{t}$ is not unambiguously
context-free.  We also show that this result holds for
generalized Thue--Morse words as well.

For an infinite word $\mathbf{w}$, let $p_{\mathbf{w}}(n)$ denote
the \emph{subword complexity function} of $\mathbf{w}$.  That is,
the value of $p_{\mathbf{w}}(n)$ is equal to the number of subwords
of length $n$ that occur in $\mathbf{w}$.  Let $L_{\mathbf{w}}$
denote the set of words over the alphabet of $\mathbf{w}$
that are not subwords of $\mathbf{w}$.

Brlek \cite{Brl89} and de Luca and Varricchio \cite{LV89} (see
also the subsequent work of Avgustinovich \cite{Avg94}, Tapsoba \cite{Tap94},
Frid \cite{Fri97}, and Tromp and Shallit \cite{TS95}) determined that
\begin{equation}
\label{subwords}
p_{\mathbf{t}}(n+1) =
\begin{cases}
2 & \text{if $n=0$}, \\
4 & \text{if $n=1$}, \\
4n - 2^a & \text{if $n = 2^a + b$, where $a \geq 1$, $0 \leq b < 2^{a-1}$}, \\
4n - 2^a - 2b & \text{if $n = 2^a + 2^{a-1} + b$, where $a \geq 1$,
$0 \leq b < 2^{a-1}$}.
\end{cases}
\end{equation}

Based on this characterization, we prove the following theorem.

\begin{theorem}
\label{tm}
The language $L_{\mathbf{t}}$ is not unambiguously context-free.
\end{theorem}

\begin{proof}
Let
\[
F(X) = \sum_{n \geq 1} p_{\mathbf{t}}(n) X^n
\]
be the generating series of the subwords of the Thue--Morse word.
We show that $F(X)$ is transcendental over $\mathbb{Q}(X)$.
Suppose to the contrary that $F(X)$ is algebraic.  Then the series
\[
G(X) = \sum_{n \geq 1} (p_{\mathbf{t}}(n+1) - p_{\mathbf{t}}(n)) X^n,
\]
whose coefficients form the sequence of \emph{first differences} of
$p_{\mathbf{t}}(n)$, is also algebraic.  Note that for all $n \geq 1$,
\[
p_{\mathbf{t}}(n+1) - p_{\mathbf{t}}(n) \leq 4,
\]
so that the coefficients of $G(X)$ are bounded.
Applying Fatou's Theorem to $G(X)$, we see that $G(X)$ is either
rational or transcendental.  By assumption, $G(X)$ is algebraic,
so it must be rational.  But then the sequence
\[
\Delta = (p_{\mathbf{t}}(n+1) - p_{\mathbf{t}}(n))_{n \geq 1}
\]
is ultimately periodic.  We easily verify from (\ref{subwords})
that this is not the case: for instance, $\Delta$ contains
arbitrarily large ``runs'' of $4$'s.  This contradiction implies
the transcendence of $F(X)$.

Alternatively, one may note that the series $H(X)$, whose coefficients
form the sequence of \emph{second differences} of $p_{\mathbf{t}}(n)$, is a
gap series, and one may therefore apply Hadamard's gap theorem to $H(X)$.

The desired result follows by applying the Chomsky--Sch\"utzenberger Theorem.
\end{proof}

\noindent\textbf{Note}: The use of analytic techniques in the
proof of Theorem~\ref{tm} may be avoided by applying instead
the theorems of Christol \cite{Chr79,CKMR80} and Cobham \cite{Cob69}.
See the paper of Allouche \cite{All99} for some examples of this
approach.

Next we consider generalizations of the Thue--Morse word.
Let $s_2(n)$ denote the sum of the digits in the base-$2$ expansion
of $n$.  It is well known that the Thue--Morse word
$\mathbf{t} = t(0)t(1)t(2)\cdots$ can be defined by $t(n) = s_2(n) \bmod 2$.
For $k \geq 2$, we define the \emph{generalized Thue--Morse word}
$\mathbf{t}_k$ by $\mathbf{t}_k(n) = s_2(n) \bmod k$, so that
$\mathbf{t} = \mathbf{t}_2$.  Tromp and Shallit \cite{TS95}
characterized the subword complexity of these words as follows:
\[
p_{\mathbf{t}_k}(n+1) =
\begin{cases}
k & \text{if $n=0$}, \\
k^2 & \text{if $n=1$}, \\
k(kn - 2^{a-1}) & \text{if $n = 2^a + b$, where $a \geq 1$, $0 \leq b < 2^{a-1}$}, \\
k(kn - 2^{a-1} - b) & \text{if $n = 2^a + 2^{a-1} + b$, where $a \geq 1$,
$0 \leq b < 2^{a-1}$}.
\end{cases}
\]
One therefore proves the following result in a manner entirely
analogous to that of Theorem~\ref{tm}.

\begin{theorem}
For $k \geq 2$, the language $L_{\mathbf{t}_k}$ is not unambiguously
context-free.
\end{theorem}

\section{Discussion and future work}

To complete the work discussed here, it remains to determine whether
or not the languages $L_o(2)$, $\widetilde{L_o}(2)$, and
$L_{\mathbf{t}_k}$ are context-free.  We discuss some related issues below.

Moss\'e \cite{Mos96} and Frid \cite{Fri97,Fri97b,Fri98} have written
several papers showing that a large class of words generated by
iterating morphisms have subword complexity functions that behave
similarly to that of the Thue--Morse word; i.e., they are piecewise linear
on exponentially growing intervals.  The first difference sequence
of such subword complexity functions is therefore either constant
or not ultimately periodic.  If it were possible to characterize
those words for which the latter situation occurs, one might
generalize the argument of Theorem~\ref{tm} to a larger class
of words.

One may also apply this argument in cases where the subword complexity
function is not linear.  As an example, we may consider the generating
series of the \emph{paperfolding words}.  Allouche and Bousquet-M\'elou
\cite{AB94} have shown that the number $f(n)$ of subwords of length $n$
of the paperfolding words is given by
\[
f(n) =
\begin{cases}
2^n & \text{if $1 \leq n \leq 3$}, \\
2^{a+1}n - 5\cdot 4^{a-1} & \text{if $2^a \leq n < 2^a + 2^{a-1}$,
where $a \geq 2$}, \\
3\cdot 2^a n - 11\cdot 4^{a-1} & \text{if $2^a + 2^{a-1} \leq n < 2^a +
2^{a-1} + 2^{a-2}$, where $a \geq 2$}, \\
2^{a+1}n - 4^a & \text{if $2^a + 2^{a-1} + 2^{a-2} \leq n < 2^{a+1}$,
where $a \geq 2$},
\end{cases}
\]
and they have shown that the corresponding generating function $F(X)$ is
transcendental.  We may also deduce the transcendence of $F(X)$ by
considering the series $H(X)$, whose coefficients form the second
difference sequence of $f(n)$.  Noting that $H(X)$ is a gap series,
we may apply Hadamard's gap theorem to derive the desired result.

Noting that $F(X)$, along with the other generating functions considered
earlier, has coefficients that are polynomially bounded, we take this
opportunity to mention a remarkable recent result of D'Alessandro,
Intrigila, and Varricchio \cite{AIV06}: \emph{If a context-free
language has only polynomially many words of length $n$, then
its generating function is rational}.  Applying this result to
the paperfolding words, for instance, one recovers a result of Lehr
\cite{Leh92,AB94}, namely, that the set of subwords of the paperfolding words
is not context-free.

\section{Acknowledgments}
Thanks to Jeffrey Shallit for reading an earlier draft of this work and
offering some helpful suggestions.

\end{document}